\documentclass[11pt]{article}
\usepackage{epsfig}

\input amssymb.sty
\begin{document}
\sloppy

\font\amsy=msbm10

\newtheorem{defn}{Definition}[section]
\newtheorem{example}[defn]{Example}
\newtheorem{question}[defn]{Question}
\newtheorem{prop}[defn]{Proposition}
\newtheorem{thm}[defn]{Theorem}
\newtheorem{lem}[defn]{Lemma}
\newtheorem{cor}[defn]{Corollary}
\newtheorem{conj}[defn]{Conjecture}
\newtheorem{rem}[defn]{Remark}

\newcommand{\dickebox}{{\vrule height5pt width5pt depth0pt}}
\def\NN{\mathbb{N}}
\def\ZZ{\mathbb{Z}}
\def\QQ{\mathbb{Q}}
\def\CC{\mathbb{C}}

\font\amsy=msbm10
\def\IC{\hbox{\amsy\char'103}}
\def\sg{\sigma}
\def\LL{\Cal L}
\def\BB{\Cal B}
\def\AA{\Cal A}
\def\lm{\lambda}
\def\KK{\Cal K}
\def\dim{\mbox{\rm dim }}
\def\sgn{\mbox{\rm sgn}}
\def\om{\omega}
\def\Om{\Omega}
\def\IC{\hbox{\amsy\char'103}}
\font\diff=msbm10
\def\leneq{\hbox{\diff\char'010}}
\def\geneq{\hbox{\diff\char'011}}
\setlength{\parindent}{0pt}

{\Large \bf
\begin{center} Expected lengths and distribution functions for
Young diagrams in the hook
\end{center}}

\medskip

\centerline{June 19, 2005}
\bigskip
\centerline{A.~Regev\footnote{Partially supported by Minerva Grant
No. 8441. 2000 Mathematics Subject Classification. Primary 60C05;
Secondary 45E05, 05A05.}}

\bigskip

\centerline{The Weizmann Institute of Science, Rehovot 76100,
Israel}
\medskip
\centerline{amitai.regev@weizmann.ac.il} \vskip 0.8cm

\begin{abstract}
We consider $\beta$--Plancherel measures \cite{Ba.Ra.} on subsets
of partitions -- and their asymptotics. These subsets are the
Young diagrams contained in a $(k,\ell)$--hook, and we calculate
the asymptotics
 of the expected shape of these diagrams, relative to such
 measures. We also calculate the asymptotics of the distribution
 function of the lengths of the rows and the columns for these diagrams. This
 might be considered as the restriction to the $(k,\ell)$--hook of the 
fundamental
 work of Baik, Deift and Johansson~\cite{B.D.J.1}. The above
 asymptotics are given here by ratios of certain Selberg-type
 multi--integrals.
\end{abstract}
\medskip

\centerline{ \centerline}
\bigskip

\section{Introduction} \label{sec1}
This paper studies the asymptotics of certain
"$\beta$--Plancherel" measures on subsets of partitions. Let
$Y_n=\{\lm\mid\lm\vdash n\}$ denote the partitions of $n$, and let
$f^\lm$ denote the number of standard Young tableaux of shape
$\lm=(\lm_1,\lm_2,\ldots)$. For general references regarding
partitions, Young diagrams and Young tableaux -- see
\cite{Mac},~\cite{St1}. Let $0<\beta \in \mathbb{R}$. Baik and
Rains \cite{Ba.Ra.} consider the following "$\beta$--{\it
Plancherel}" measure $M^{\beta}_n$ on $Y_n$:
\begin{eqnarray}\label{br1}
M^{\beta}_n(\lm):=\frac{(f^\lm)^\beta}{\sum_{\mu\vdash
n}(f^\mu)^\beta}.\end{eqnarray} Indeed, $M^2_n$ is the so called
Plancherel measure on $Y_n$. We generalize to subsets
$\Gamma_n\subseteq Y_n$, considering the subset of the partitions
(i.e. diagrams) in the $(k,\ell)$ hook: Let $k,\,\ell\ge 0$ be
integers and let $\Gamma_n=H(k,\ell;n)$ denote the following
subset of $Y_n$:
\[H(k,\ell;n)=\{\lm\vdash n\mid \lm_{k +1}\le\ell\}.\]
These subsets arise in the representation theory of Lie groups,
algebras and superalgebras, see for example \cite{Be.Re.}. The
measures $\rho ^{(\beta;k,\ell)}_n$ below are the $(k,\ell)$--hook
restrictions of the above measures $M_n^\beta$.
\begin{defn}\label{st11}
Let $\lm\in H(k,\ell;n)$ and $\beta>0$, then
\[\rho ^{(\beta;k,\ell)}(\lm)=\rho ^{(\beta;k,\ell)}_n(\lm):=
\frac{(f^\lm)^\beta}{\sum_{\mu\in H(k,\ell;n)}(f^\mu)^\beta}.
\]
\end{defn}

\medskip

\subsection{Expected shape}
Given $\Gamma_n\subseteq Y_n,~~n=1,2,\ldots$ and the probability
measures $\rho=\{\rho_n\}_{n=1}^\infty$ on the $\Gamma_n$--s, one
studies the asymptotics of the expected value (i.e.~average
length) of the first row $\lm_1$, denoted $\lm_{1,E}$, and
similarly for the second row $\lm_{2,E}$, etc. Similarly for the
columns. Explicitly, when $\Gamma_n=H(k,\ell;n)$, expected values
are given by the following definition.

\begin{defn}\label{def1} If $\lm\vdash n$, we write
$\lm=(\lm_{1,n},\lm_{2,n},\ldots)$. Also, $\lm'$ is the conjugate
partition of $\lm$. Let $1\le p\le k$ and $1\le q\le \ell$. The
expected value of the $p$-th row is $E(\lm_p)=\lm^{(\beta;k,\ell)}
_{p,E}(n)$, where
\[\lm^{(\beta;k,\ell)} _{p,E}(n)=\frac{\sum_{\lm\in
H(k,\ell;n)}\lm_{p,n} \cdot(f^\lm)^{\beta}} {\sum_{\lm\in
H(k,\ell;n)}(f^\lm)^{\beta}}\quad\mbox{and}\quad
  \lm^{(\beta;k,\ell)}_{E}(n)=(\lm^{(\beta;k,\ell)}_{1,E}(n),\lm^{(\beta;k,\ell)}_{2,E}(n),\ldots)\]

Similarly for the expected $q$-th column
\[\lm'^{(\beta;k,\ell)} _{q,E}(n)=\frac{\sum_{\lm\in
H(k,\ell;n)}\lm'_{q,n} \cdot(f^\lm)^{\beta}} {\sum_{\lm\in
H(k,\ell;n)}(f^\lm)^{\beta}} \quad\mbox{and}\quad
\lm^{'(\beta;k,\ell)}_{E}(n)=(\lm^{'(\beta;k,\ell)}_{1,E}(n),\lm^{'(\beta;k,\ell)}_{2,E}(n),\ldots)
.\]
\end{defn}

Of course, one can replace $H(k,\ell;n)$ in the above definition
by other subsets $\Gamma_n\subseteq Y_n$.

\bigskip

The  case $\Gamma_n=Y_n$ and $\beta=2$ (Plancherel) has a long
history. Let
\[w(n)=\frac{\sum_{\lm\vdash n}
\lm_{1,n} \cdot(f^\lm)^{2}} { \sum_{\lm\vdash n} (f^\lm)^{2}
}=\frac{\sum_{\lm\vdash n} \lm_{1,n} \cdot(f^\lm)^{2}} {n!} \] be
the expected value of the first row -- for the Plancherel measure
$M^2_n$. Hammersley \cite{H} showed that the limit
$c=\lim_{n\to\infty} w(n)/\sqrt n~$ exists. Vershik and Kerov
\cite{VK1} proved that $c=2$ (independently, Logan and Shepp
\cite{LS} proved that $c\le2$). Vershik and Kerov -- and Logan and
Shepp also determined the asymptotics of the expected shape $\lm$
in this case.

\medskip

Recently, in a major breakthrough paper, Baik, Deift and Johansson
\cite{B.D.J.1} determined the distribution function of the
asymptotics of the first row, relating it to the {\it Tracy-Widom
distribution} \cite{TW1}, see also \cite{Ba.Ra.} and \cite{TW2}.
The distribution function for the second row is given, by these
same authors, in \cite{B.D.J.2}. The distribution functions for
the general rows are given in \cite{B.O.O.}, \cite{Jo} and
\cite{Ok}; see also \cite{Bo} for the analogue results for colored
permutations. The above results also establish deep connections
with the theory of {\it random matrices} \cite{Mehta}. For
detailed reviews of these results -- see \cite{A.D.}
and~\cite{St2}.

\medskip

\bigskip

~~~The main objective of the present paper is to compute the
asymptotics of the above expected values (i.e.~shapes)
$\lm^{(\beta;k,\ell)}_{E}(n)$, as well as the corresponding
distribution functions. The first term approximation is relatively
simple, as we show that for each $1\le p\le k$ and $1\le q\le
\ell$
\[\lm^{(\beta;k,\ell)}_{p,E}(n),~\lm'^{(\beta;k,\ell)}_{q,E}(n)\sim\frac{n}{k+\ell},\] see
Theorem~\ref{first}. Second term approximations of
$\lm^{(\beta;k,\ell)}_{p,E}(n)$ are introduced and studied in
Sections~\ref{2d1}, \ref{2d2}, and they have different values for
different rows and for different columns. These
second-term-approximations are given as ratios of certain
Selberg-type integrals, see Theorems~\ref{ratio} and~\ref{ratio2}
below.

\subsection{Distribution functions}\label{6}
In Section \ref{sec6} we introduce and study the asymptotics of
the Distribution functions
$\lm^{(\beta;k,\ell)}_{p}(n,z),\;\lm'^{(\beta;k,\ell)}_{q}(n,z)$
for the lengths of the rows and the columns in $H(k,\ell;n)$ --
with respect to the above measures. We are able to calculate,
asymptotically, a first--term approximation of these functions,
but only a conjecture is given, about the second term
approximations. That first--term approximation is
\[\lm^{(\beta;k,\ell)}_{p}(n,z),\;\lm'^{(\beta;k,\ell)}_{q}(n,z)
\sim\frac{n}{k+\ell}\cdot r_{(k,\ell),\beta}(z),\] where
$r_{(k,\ell),\beta}(z)$ is given by Equation~(\ref{dist222}), see
Theorem~\ref{first222}.

\subsection{Comparison with maximal shape}\label{1.2}
 Given a subset of partitions $\Gamma_n\subseteq Y_n$, one looks for
$\lm\in\Gamma_n$ with maximal degree $f^\lm$. Call it {\it maximal
shape} (with respect to $\Gamma_n$) and denote it by $\lm_{max}$.
In Sections~\ref{sec8},~\ref{sec9} and~\ref{sec10} the expected
shapes for $\beta =1,2$ are compared with the maximal shape. When
$\Gamma_n=Y_n$, the asymptotics of $\lm_{max}$ was calculated by
Vershik and Kerov \cite{VK1}, \cite{VK2},  and by Logan and Shepp
\cite{LS}. In particular, they proved that asymptotically, the
{\it expected} shape for $\beta=2$ and the {\it maximal} shape are
the same, and that shape is given by the two axes and by the curve
\begin{eqnarray}\label{gmax} y=1+\left(\frac{2}{\pi}
\right)[x\cdot\sqrt{1-x^2}- arccos\,x].
\end{eqnarray}

\medskip

A comparison with the case $\Gamma_n=H(k,\ell;n)$ is intriguing.

\medskip

When $\Gamma_n$ is the $k$--strip $\Gamma_n=H(k,0;n)$, the
asymptotics of $\lm$ with maximal $f^\lm$ was calculated in
\cite{AR}, and is given by the curve
\begin{eqnarray}\label{smax} y=\left(\frac{2}{\pi}
\right)[x\cdot arcsin\, x+\sqrt{1-x^2}].
\end{eqnarray} When $\Gamma_n=H(k,\ell;n)$, the maximal $\lm$ was given in
\cite{hook}. These results are reviewed in Section~\ref{max5}.
Consider for example the 'strip' case $\ell=0$, and denote the
maximal $\lm$ by $\lm^{(k,0)}_{max}$. Comparing it with
$\lm^{(2;k,0)}_{E}$, the 'Plancherele' expected $\lm$ in
$H(k,0;n)$, we show that these asymptotic shapes are {\it not}
equal -- even in their first raw. Nevertheless, {\it numerically}
$\lm^{(k,0)}_{max}$ and $\lm^{(2;k,0)}_{E}$ are remarkably close,
at least in the few special cases we check below, see
Section~\ref{sec8} .

\medskip

Also, the asymptotics of $\lm_{max}$ for $\Gamma_n=Y_n$ is {\it
not} the limit case of $\lm^{(k,0)}_{max}$ as $~k\to\infty$, but
the similarity between (\ref{gmax}) and (\ref{smax}) is
intriguing. It should be interesting to see if the ratios of the
Selberg-type integrals, which give the expected shapes and the
distribution functions for $\Gamma_n=H(k,\ell;n)$, are in any way
related to the Tracy--Widom distributions \cite {TW1} \cite{TW2},
which give the (Plancherel and the 'involution') distribution
functions in $Y_n$.

\subsection{RSK}\label{rsk}
In the case of $\rho ^{(1;k,\ell)}_n$ and $\rho ^{(2;k,\ell)}_n$,
the RSK correspondence provides an interesting interpretation of
the above asymptotics. The RSK (Robinson--Schensted--Knuth)
correspondence $\sg\longleftrightarrow(P_\lm,Q_\lm)$ corresponds
$\sg\in S_n$ with a pair of standard Young tableaux of shape
$\lm$~~\cite{St1}. In the Plancherel case $\beta=2$ it relates the
above expected values of the first row to the statistics of the
longest increasing (and decreasing) subsequences in permutation.
For example, when $\sg\longleftrightarrow(P_\lm,Q_\lm)$, $\lm_1$
is the length of a longest increasing subsequence in $\sg$, while
$\lm'_1$ is the length of a longest decreasing subsequence in
$\sg$. By C. Green's theorem there are similar interpretations for
$\lm_2,\;\lm_3,$ etc. For a detailed account of the RSK see
\cite{St1}. Thus the results in \cite{B.D.J.1} etc.~can also be
stated in terms of longest increasing subsequences in
permutations.

\medskip

It is well known that $\sg$ is an involution if''f
$\sg\longleftrightarrow (P_\lm,P_\lm)$. The analogue Probability
theory of longest increasing subsequences in involutions in $S_n$
is done in \cite{Ba.Ra.}.

\medskip

 Denote by
$S_{k,\ell;n}\subseteq S_n$ the subset of the permutations $\sg\in
S_n$ such that under the RSK correspondence
$\sg\longleftrightarrow (P_\lm ,Q_\lm)$, we have $\lm\in
H(k,\ell;n)$.
For example, $S_{k,0;n}$ is the subset of those permutations in
$S_n$ where any descending subsequence has length $\le k$.
Thus, $\lm_{1,E}^{(1;k,\ell)}$ is the expected value of the
longest increasing subsequence in the involutions in
$S_{k,\ell;n}$.


\medskip

\section{Selberg type integrals}\label{Selberg} As mentioned above, the main
results in this paper involve Selberg--type integrals, hence we
briefly review these type of multi--integrals. In~\cite{Se} A.
Selberg proved the following formula:
\[\int_0^1\cdots\int_0^1(u_1\cdots u_n)^{x-1}\cdot [(1-u_1)\cdots
(1-u_n)]^{y-1}\cdot\prod_{1\le i<j\le n}|u_i-u_j|^{2z}du_1\cdots
du_n=~~~~\]
\[~~~~~~~~~~~~~~~~~~~~~~~~~~~~~~~~~~~~=\prod_{k=1}^n
\frac{\Gamma (1+kz)\cdot\Gamma (x+(k-1)z)\cdot \Gamma
(y+(k-1)z)}{\Gamma (1+z)\cdot\Gamma (x+y+(n+k-2)z)}.\]

Various integral formulas can be deduced from Selber's integra,
see for example~\cite{Mac2} for the Macdonald--Mehta integras. For
example Mehta's integral formula (which was a conjecture for some
time)
\[ \int_{\mathbb{R}^k}e^{-(1/2)(\sum x_i^2)}\prod_{1\le i<j\le k}
|x_i-x_j|^{2z}dx_1\cdots dx_k=(\sqrt
{2\pi})^k\cdot\prod_{j=1}^k\frac{\Gamma(1+jz)}{\Gamma(1+z)}
\]
can be deduced from Selberg's formula, see~\cite{Mehta} for
details. We call these and related integrals "Selberg--type
integrals". A connection between the RSK and these integrals, as
well as with {\it random matrices}, appears in~\cite{Strip},
\cite{Be.Re.}. Since the formulas from~\cite{Strip}, \cite{Be.Re.}
are needed later, we record it here, together with certain
variations of these asymptotics and integrals that are also needed
below. Let
\[\Omega_k=\{(x_1,\ldots,x_k)\in \mathbb{R}^k\mid x_1\ge
x_2\ge\cdots \ge x_k ~\mbox{and}~x_1+\cdots +x_k=0\},\] and more
generally,
\[\Omega_{(k,\ell)}=\{(x_1,\ldots,x_k,y_1,\ldots,y_\ell)\mid
x_1\ge\cdots\ge x_k\,;~~ y_1\ge\cdots\ge y_\ell\,;~~\sum x_i+\sum
y_j=0\}.\]

\begin{thm}\label{sheli1}(Theorem 2.10 in~\cite{Strip})
Let $\gamma_k=(1/\sqrt{2\pi})^{k-1}\cdot k^{k^2/2}$ and
$D_k(x)=\prod_{1\le i<j\le k}(x_i-x_j)$, then
\[\sum_{\lm\in H(k,0;n)}(f^\lm)^\beta\sim\left [\gamma_k\cdot
\left (\frac{1}{n}\right )^{(k-1)(k+2)/4}\cdot k^n  \right ]^\beta
\cdot(\sqrt n)^{k-1}\cdot I(k,0,\beta),\] where
\[ I(k,0,\beta)=\int_{\Omega_k}\left [D_k(x)\cdot e^{-\frac{k}{2}(\sum
x_i^2)}
  \right ]^\beta~d^{(k-1)}x.\]
\end{thm}
Note that by certain symmetry  properties of the above integrand,
$\Omega_k$ is transformed in \cite{Strip} into $\mathbb{R}^k $.

\medskip

Given $\lm=(\lm_{1,n},\lm_{2,n},\ldots )\vdash n$, write
$\lm_{p,n}=n/k+c_{p,n}\cdot\sqrt n$ and denote
$c_{p,n}=c_{p,n}(\lm)$. The same arguments in \cite{Strip} prove
the following theorem.
\begin{thm}\label{sheli2}
\[\sum_{\lm\in H(k,0;n)}c_{p,n}(\lm)\cdot(f^\lm)^\beta\sim\left [\gamma_k\cdot
\left (\frac{1}{n}\right )^{(k-1)(k+2)/4}\cdot k^n  \right ]^\beta
\cdot(\sqrt n)^{k-1}\cdot I^*(k,0,\beta),\] where
\[ I^*(k,0,\beta)=\int_{\Omega_k} x_p\cdot \left [D_k(x)\cdot e^{-\frac{k}{2}(\sum
x_i^2)}
  \right ]^\beta~d^{(k-1)}x.\]
\end{thm}

The $(k,\ell)$--hook analogue of Theorem~\ref{sheli1} is proved
in~\cite{Be.Re.}:
\begin{thm}\label{sheli3}(See Theorem 7.18 in~\cite{Be.Re.})
Let \[\gamma_{k,\ell}=(1/\sqrt{2\pi})^{k+\ell-1}\cdot
(k+\ell)^{(k^2+\ell^2)/2}\cdot(1/2)^{k\ell},\] then
\[\sum_{\lm\in H(k,\ell;n)}(f^\lm)^\beta\sim
\left [\gamma_{k,\ell}\cdot\left (\frac{1}{n}\right
)^{(k(k+1)+\ell(\ell+1)-2)/4} \cdot (k+\ell)^n \right] ^\beta\cdot
(\sqrt n)^{k+\ell-1}\cdot I(k,\ell,\beta),\] where
\[I(k,\ell,\beta)=\int_{\Omega_{k,\ell}}\left [D_k(x)\cdot D_\ell (y)\cdot
e^{-\frac{k+\ell}{2}(\sum x_i^2+\sum y_j^2)} \right ]^\beta
~d^{(k+\ell-1)}(x;y).\]
\end{thm}
The same arguments also prove
\begin{thm}\label{sheli4}
Let $1\le p\le k$, let $\lm=(\lm_{1,n},\lm_{2,n},\ldots )\in
H(k,\ell;n)$ and define $c_{p,n}(\lm)$
\[via: \quad\lm_{p,n}=\frac{n}{k+\ell}+c_{p,n}(\lm)\cdot\sqrt n.\quad\mbox {Similarly
for}\quad 1\le q\le \ell\quad\mbox{and}\quad
c'_{q,n}(\lm):=c_{q,n}(\lm').\] Then
\[\sum_{\lm\in H(k,\ell;n)}c_{p,n}(\lm)\cdot
(f^\lm)^\beta\sim~~~~~~~~~~~~~~~~~~~~~~~~~~~~~~~~~~~~~~~~~~~~~~~~~~~~~~~~~~~~~~~~~~~~~~~~\]
\[~~~~\left [\gamma_{k,\ell}\cdot\left (\frac{1}{n}\right
)^{(k(k+1)+\ell(\ell+1)-2)/4} \cdot (k+\ell)^n \right] ^\beta\cdot
(\sqrt n)^{k+\ell-1}\cdot I^*(k,\ell,\beta),\] where
\[I^*(k,\ell,\beta)=\int_{\Omega_{k,\ell}}x_p\cdot \left [D_k(x)\cdot D_\ell (y)\cdot
e^{-\frac{k+\ell}{2}(\sum x_i^2+\sum y_j^2)} \right ]^\beta
~d^{(k+\ell-1)}(x;y).\] Similarly for the sum \[\sum_{\lm\in
H(k,\ell;n)}c'_{q,n}(\lm)\cdot (f^\lm)^\beta,\] with the
corresponding integral
\[I^{'\,*}(k,\ell,\beta)=\int_{\Omega_{k,\ell}}y_q\cdot \left
[D_k(x)\cdot D_\ell (y)\cdot e^{-\frac{k+\ell}{2}(\sum x_i^2+\sum
y_j^2)} \right ]^\beta ~d^{(k+\ell-1)}(x;y).\]
\end{thm}
Generalizations of Theorem~\ref{sheli4} -- with more general
functions of the $c_{p,n}(\lm),\;c'_{q,n}(\lm)$ -- are rather
obvious, but will not be given here.

\medskip


The same arguments of Section 7 of~\cite{Be.Re.}, applied to the
asymptotics of both numerator and denominator, prove the following
theorem -- which is needed later.
\begin{thm}\label{sheli5}
Let $z>0$, \[H(k,\ell;n,z)=\{\lm\in H(k,\ell;n)\mid
 \lm_{1,n},\;\lm_{1,n}\le\frac{n}{k+\ell}+z\sqrt n\,\}\] and let
\[\Omega_{(k,\ell),z}=\{(x_1,\ldots,x_k;y_1,\ldots,y_\ell)\in\Omega_{(k,\ell)}\mid x_1,\;y_1\le z\}.\]
Then \[\frac{\sum_{\lm\in H(k,\ell;n,z)
 }(f^\lm)^{\beta}} {\sum_{\lm\in H(k,\ell;n)
}(f^\lm)^{\beta}}\sim
~~~~~~~~~~~~~~~~~~~~~~~~~~~~~~~~~~~~~~~~~~~~~~~\]\[~~~~~~~~~~~~~~~~~~~~\sim\frac{\int_{\Omega_{(k,\ell),z}}
\left [ D_k(x)\cdot D_\ell(y)\cdot e^{-\frac{k+\ell}{2}(\sum
x_i^2+\sum y_j^2)}\right ]^\beta
d^{(k+\ell-1)}(x;y)}{\int_{\Omega_{(k,\ell)}} \left [D_k(x)\cdot
D_\ell(y) \cdot e^{-\frac{k+\ell}{2}(\sum x_i^2+\sum y_j^2)}\right
]^\beta d^{(k+\ell-1)}(x;y)}.\]
\end{thm}

We call all the above   "{\it Selberg--type integrals}", and
remark that the above expected values and distribution functions
are given below as ratio of such integrals.

\section{The main results}
\subsection{The expected values}
We study the expected values of the row and of the column lengths
in $\Gamma_n=H(k,\ell;n)$ with respect to the measures
$\rho_n^{(\beta;k,\ell)}$ introduced in Definition~\ref{st11}. The
first term asymptotics is given by
\begin{thm} (See Theorem \ref{first})
Let $\beta > 0$ and $\Gamma_n=H(k,\ell;n)$. For each $1\le p\le
k$,
\[\lm^{(\beta;k,\ell)}_{p,E}(n)\sim\frac{n}{k+\ell},\quad\mbox{
namely}\quad \lim_{n\to\infty}\left
({\lm^{(\beta;k,\ell)}_{p,E}(n)}\right )/{ \left( \frac{n}{k+\ell}
\right)}=1.\] Similarly, for each $1\le q\le \ell$,
\[\lm'^{(\beta;k,\ell)}_{q,E}(n)\sim\frac{n}{k+\ell}.\]
\end{thm}

~~The second term approximations are given as follows. Define
$c^{(\beta;k,\ell)}_{p,E}(n)$ and $c^{(\beta;k,\ell)}_{p,E}$ via
\[\lm^{(\beta;k,\ell)}_{p,E}(n)=\frac{n}{k+\ell}+c^{(\beta;k,\ell)}_{p,E}(n)\cdot\sqrt
n,\quad\mbox{and}\quad c^{(\beta;k,\ell)}_{p,E}=\lim
_{n\to\infty}c^{(\beta;k,\ell)}_{p,E}(n).\] Similarly for the
columns. Then
\begin{thm} (see Theorem~\ref{ratio2})
 Let $1\le p\le k$,  then the limit $c^{(\beta;k,\ell)}_{p,E}$  exists,
 and is given as follows:
\begin{eqnarray}\label{hookeqn111}
c^{(\beta;k,\ell)}_{p,E}=\frac{\int_{\Omega_{k}} x_p\cdot \left
[D_k(x_1,\ldots,x_k) \cdot e^{-\frac{k+\ell}{2}(x_1^2+\cdots
+x_k^2)}\right ]^\beta d^{(k-1)}(x)}{\int_{\Omega_{k}} \left
[D_k(x_1,\ldots,x_k) \cdot e^{-\frac{k+\ell}{2}(x_1^2+\cdots
+x_k^2)}\right ]^\beta d^{(k-1)}(x)}.
\end{eqnarray}
 Similarly for the  columns.
\end{thm}
 Equation~(\ref{hookeqn111}) (or~(\ref{hookeqn1})) is a consequence of the seemingly more symmetric
 Equation~(\ref{eqnhook1}).

\medskip

In Sections~\ref{sec8},~\ref{sec9} and~\ref{sec10} the expected
shapes for $\beta =1,2$ are compared with the maximal shape
$\lm_{max}$ in few special cases. As mentioned in
Section~\ref{1.2}, this shows a different behavior in the hook
case $\Gamma_n=H(k,\ell;n)$ compared with the general case
$\Gamma_n=Y_n$.

\subsection{The distribution function for the first row}
In Section~\ref{sec6} we study the distribution functions of the
length on the rows and the columns. We calculate the first term
approximations and conjecture the second term approximations.
Given $0<z\in\mathbb{R}$, denote
\[H(k,\ell;n,z)=\{\lm\in H(k,\ell;n)\mid
 \lm_{1,n},\;\lm'_{1,n}\le\frac{n}{k+\ell}+z\sqrt n\,\}.\]
 Let $1\le p\le k,~1\le q\le\ell$.
The distribution of the length of the $p$--th row as a function of
$z$ is defined  as

\[\lm^{(\beta;k,\ell)} _p(n,z)=\frac{\sum_{\lm\in H(k,\ell;n,z)
 }\lm_{p,n}
\cdot(f^\lm)^{\beta}} {\sum_{\lm\in H(k,\ell;n)
}(f^\lm)^{\beta}},\] and similarly for the columns. Recall
$\Omega_{(k,\ell)}$ from Section~\ref{Selberg} and denote
\[\Omega_{(k,\ell),z}=\{(x_1,\ldots,x_k;y_1,\ldots,y_\ell)\in\Omega_{(k,\ell)}\mid x_1,\;y_1\le z\}.\]
Then
\begin{thm}(see Theorem~\ref{first222})~~  Let $\beta,\,z > 0$ and
denote
\begin{eqnarray}  r_{(k,\ell),\beta}(z)=
\frac{\int_{\Omega_{(k,\ell),z}} \left [ D_k(x)\cdot
D_\ell(y)\cdot e^{-\frac{k+\ell}{2}(\sum x_i^2+\sum y_j^2)}\right
]^\beta d^{(k+\ell-1)}(x;y)}{\int_{\Omega_{(k,\ell)}} \left
[D_k(x)\cdot  D_\ell(y) \cdot e^{-\frac{k+\ell}{2}(\sum x_i^2+\sum
y_j^2)}\right ]^\beta d^{(k+\ell-1)}(x;y)}.\end{eqnarray}   Then
\[\lm^{(\beta;k,\ell)}_{p}(n,z),\;\lm'^{(\beta;k,\ell)}_{q}(n,z)
\sim\frac{n}{k+\ell}\cdot r_{(k,\ell),\beta}(z).\]
\end{thm}
In Section~\ref{sec7} we make some conjectures about the second
term approximations of $\lm^{(\beta;k,\ell)}_{p}(n,z)$ and
$\lm'^{(\beta;k,\ell)}_{q}(n,z)$, both in terms of ratios of
Selberg--type integrals. In the last three sections (Sections 7, 8
and 9), we calculate some special cases and include also some
computer calculations.

\bigskip


\part{Expected shape, the $\beta$--Plancherel probability}
\section{First term approximation} Recall the notation $\sim$:
Let $a_n,\;b_n$ be two sequences of, say, real numbers, and assume
$b_n\ne 0$ if $n$ is large enough. Then $a_n\sim b_n$ if
$\lim_{n\to\infty}a_n/b_n=1$. Extend $\sim$ to vectors as follows:
$~(a_{1,n},\,\ldots, a_{r,n})\sim (b_{1,n},\,\ldots, b_{r,n})~$
if''f $~a_{i,n}\sim b_{i,n}$ for $i=1,\ldots r$.

\medskip

Following the techniques and arguments in Section 7 of
\cite{Be.Re.}, (see also ~\cite{Strip}) we show below the
following first approximation of the expected shape
$\lm_{E}=\lm^{(\beta;k,\ell)}_{E}$.

\begin{thm}\label{first}
Let $\beta > 0$ and $\Gamma_n=H(k,\ell;n)$ and let
$\lm^{(\beta;k,\ell)}_{p,E}(n)$ and
$\lm'^{(\beta;k,\ell)}_{q,E}(n)$ be given by Definition
\ref{def1}. For each $1\le p\le k$,
\[\lm^{(\beta;k,\ell)}_{p,E}(n)\sim\frac{n}{k+\ell},\quad\mbox{
namely}\quad \lim_{n\to\infty}\left
({\lm^{(\beta;k,\ell)}_{p,E}(n)}\right )/{ \left( \frac{n}{k+\ell}
\right)}=1.\] Similarly, for each $1\le q\le \ell$,
\[\lm'^{(\beta;k,\ell)}_{q,E}(n)\sim\frac{n}{k+\ell}.\]
\end{thm}
{\bf Proof}. We sketch the proof for the expected row length in
the case $\ell=0$, thus showing that
\[\lm^{(\beta;k,0)}_{p,E}(n)\sim\frac{n}{k}.\]
The main point is that since $\beta>0$, both sums in the numerator
and the denominator of definition~\ref{def1} are dominated by the
summands corresponding to the partitions $\lm$, such that
$\lm_i=\frac{n}{k}+c_i\cdot\sqrt n~$ and with $c_i$--s in a
bounded interval. In other words, let $a>0$ and denote\\
\[H_a(k,0;n)=\{\lm\in H(k,0;n)\mid \lm_i=\frac{n}{k}+c_i\cdot\sqrt
n,~\mbox{where}~|c_i|\le a,~~i=1,2,\ldots, k\}.\] Then
\[\lim_{n\to\infty}\lm^{(\beta;k,0)}_{p,E}(n)=\lim_{n\to\infty}\frac{\sum_{\lm\in
H(k,0;n)}\lm_{p,n} \cdot(f^\lm)^{\beta}} {\sum_{\lm\in
H(k,0;n)}(f^\lm)^{\beta}} =\lim_{a\to\infty}\left
[\lim_{n\to\infty}\frac{\sum_{\lm\in H_a(k,0;n)}\lm_{p,n}
\cdot(f^\lm)^{\beta}} {\sum_{\lm\in
H_a(k,0;n)}(f^\lm)^{\beta}}\right].\] Writing
$\lm_{p,n}=\frac{n}{k}+c_{p,n}\cdot\sqrt n$, the expression in the
brackets equals \[\frac{\sum_{\lm\in
H_a(k,0;n)}(\frac{n}{k})(f^\lm)^{\beta}}{\sum_{\lm\in
H_a(k,0;n)}(f^\lm)^{\beta}}+\frac{\sum_{\lm\in
H_a(k,0;n)}c_{p,n}\sqrt n (f^\lm)^{\beta}}{\sum_{\lm\in
H_a(k,0;n)}(f^\lm)^{\beta}}.\] The first summand equals
$\frac{n}{k}$ while the absolute value of the second summand is
bounded by $a\sqrt n$ since all $|c_p(n)|\le a$. Since
$\frac{n}{k}+a\sqrt n\sim\frac{n}{k}$, it follows that
\[\lm^{(\beta;k,0)}_{p,E}(n)\sim\frac{n}{k},\] which completes the proof in the case
$\ell=0$.~~~~~~~~~~~~~~~~~~q.e.d.

\section{Second term approximation, the 'strip' case
($\ell=0$)}\label{2d1}

Because of Theorem~\ref{first}, we look for a more subtle
approximation of $\lm^{(\beta;k,\ell)}_{p,E}(n)$, namely, we look
for the expected deviation -- of the form $c\cdot\sqrt n$ -- from
$\frac{n}{k+\ell}$. This leads us to introduce the asymptotic
expected value $c^{(\beta;k,\ell)}_{p,E}$ below. We begin with the
'strip' case $\ell=0$. The general $(k,\ell)$--hook case is given
in the next section.

\begin{defn}\label{def3}
Let $\Gamma_n=H(k,0;n)$ and let $1\le p\le k$, with
$\lm^{(\beta;k,0)}_{p,E}(n)$  given by Definition~\ref{def1}.
Define $c^{(\beta;k,0)}_{p,E}(n)$ via the equation
 \[\lm^{(\beta;k,0)}_{p,E}(n)=\frac{n}{k}+{ c^{(\beta;k,0)}_{p,E}(n)}\cdot\sqrt
 n,\qquad\mbox
 {and} \qquad c^{(\beta;k,0)}_{p,E}=\lim _{n\to\infty} c_{p,E}^{(\beta;k,0)}(n).\]
Thus, when $n$ goes to infinity,
\[\lm^{(\beta;k,0)}_{p,E}(n)\sim\frac{n}{k}+ c^{(\beta;k,0)}_{p,E}\cdot \sqrt n.\]
\end{defn}
\begin{rem}

It is not obvious that the limit $c^{(\beta;k,0)}_{p,E}=\lim
_{n\to\infty} c^{(\beta;k,0)}_{p,E}(n)$ exists. However,
Theorem~\ref{ratio} asserts that in fact, this limit does exist.
\end{rem}

Our aim is to calculate $ c^{(\beta;k,0)}_{p,E}$, thus calculating
''the second term'' in the approximation of the expected value of
the $p$-th row--length $ \lm^{(\beta;k,0)}_{p,E}$.
Definitions~\ref{def1} and~\ref{def3} obviously  imply the
equation
%
\begin{eqnarray}\label{propdef1}
c^{(\beta;k,0)}_{p,E}(n)=\left (\frac{\sum_{\lm\in
H(k,0;n)}\lm_{p,n} \left (f^\lm\right)^\beta} {\sum_{\lm\in
H(k,0;n)}\left (f^\lm\right)^\beta}-\frac{n}{k}\right
)\frac{1}{\sqrt n} .
\end{eqnarray}
If $k=1$, Equation~(\ref{propdef1}) implies that for any $\beta
>0$, $~~c^{(\beta;1,0)}_{1,E}=0$.\\

\medskip


\begin{thm}\label{ratio}
Let $k\ge 2$.  The limit $c^{(\beta;k,0)}_{p,E}=\lim _{n\to\infty}
c^{(\beta;k,0)}_{p,E}(n)$ exists, and is given by

\begin{eqnarray}\label{stripeqn1}
~~c^{(\beta;k,0)}_{p,E}=\frac{\int_{\Omega_k} x_p\cdot \left
[D_k(x_1,\ldots,x_k)\cdot e^{-\frac{k}{2}(x_1^2+\cdots
+x_k^2)}\right ]^\beta d^{(k-1)}x}{\int_{\Omega_k} \left
[D_k(x_1,\ldots,x_k)\cdot e^{-\frac{k}{2}(x_1^2+\cdots
+x_k^2)}\right ]^\beta d^{(k-1)}x}.
\end{eqnarray}

\end{thm}
{\bf Proof}. In Equation~(\ref{propdef1}) write
$\lm_{i,n}=\frac{n}{k}+c_{i,n}\cdot\sqrt n$, and consider $\lm$--s
with $c_i$ bounded in some interval: $|c_i|\le a$ for some $a>0$.
By an argument similar to the proof of Theorem \ref{first}, it
follows that as $n\to\infty$, $~~c^{(\beta;k,0)}_{p,E}(n)$ ~~is
approximated by the ratio
\begin{eqnarray}\label{ratio4} \frac{\sum_{\lm\in H_a(k,0;n)}c_{p,n} \left
(f^\lm\right)^\beta} {\sum_{\lm\in H_a(k,0;n)}\left
(f^\lm\right)^\beta}.\end{eqnarray}
The proof now follows by applying Theorem~\ref{sheli1} to the
denominator, Theorem~\ref{sheli2} to the numerator, then
cancelling equal terms. ~~~~~~~~~~~~~~~q.e.d.

\begin{rem}\label{value1}
Note that the denominator of (\ref{stripeqn1}) is actually a
''Selberg''--or a Macdonald--Mehta -- integral, which can be
evaluated for any $\beta$. For example, let $\beta=2$. By
comparing (F.2.10) with (F.4.5.2) of ~\cite{Strip}, deduce that
\[\int_{\Omega_k} \left[ D_k(x_1,\ldots,x_k)\cdot e^{-\frac{k}{2}(x_1^2+\cdots+x_k^2)}\right ]^2
d^{(k-1)}x=~~~~~~~~~~~~~~~~~~~~~~~~~~~~~~~~~~~~~~~~~~~~~~~~~~~~~~~~~~~~~\]
\[~~~~~~~~~~~~~~~~~~~~~~~~~~~~~~~~~=\left (\sqrt{2\pi}\right )^{k-1}\cdot\left
(\frac{1}{\sqrt 2}\right )^{k^2-1}\cdot\left (\frac{1}{\sqrt
k}\right )^{k^2}\cdot 1!\cdot 2!\cdots (k-1)! .\] Similarly, by
comparing (F.2.10) with (F.4.5.1) of ~\cite{Strip}, deduce that
when $\beta=1$,
\[\int_{\Omega_k} \left[ D_k(x_1,\ldots,x_k)\cdot
e^{-\frac{k}{2}(x_1^2+\cdots+x_k^2)}\right ]
d^{(k-1)}x=~~~~~~~~~~~~~~~~~~~~~~~~~~~~~~~~~~~~~~~\]

\[~~~~~~~~~~~~~~~~~~~~~~~~~~=\left(\frac{1}{\sqrt k}\right
)^{k(k+1)}\cdot\frac{1}{k!}\cdot \left(\sqrt 2\right)
^{3k-1}\cdot\frac{1}{\sqrt\pi}\cdot\prod_{j=1}^k\Gamma
\left(1+\frac{1}{2}j\right).\] For explicit values, recall that
$\Gamma\left(\frac{3}{2}\right)=\frac{\pi}{2}$, that $\Gamma
(1)=1$, and that $\Gamma(z+1)=z\Gamma(z)$.
\end{rem}

\section{Second term approximation, the $(k,\ell)$--hook
case}\label{2d2} We turn now to the general $(k,\ell)$--hook case.
\begin{defn}\label{def??}
Let $\beta >0$ and let $1\le p\le k$ and $1\le q\le \ell$, with
$\lm^{(\beta;k,\ell)}_{p,E}(n)$ given by Definition \ref{def1}.
Define $c^{(\beta;k,\ell)}_{p,E}$ via the equation
\[\lm^{(\beta;k,\ell)}_{p,E}(n)=\frac{n}{k+\ell}+c^{(\beta;k,\ell)}_{p,E}(n)\cdot\sqrt n,\quad\mbox{and}\quad
c^{(\beta;k,\ell)}_{p,E}=\lim
_{n\to\infty}c^{(\beta;k,\ell)}_{p,E}(n).\]

Similarly for the columns:
\[\lm'^{(\beta;k,\ell)}_{q,E}(n)=\frac{n}{k+\ell}+c'^{(\beta;k,\ell)}_{q,E}(n)\cdot\sqrt n, \quad\mbox{and}\quad
c'^{(\beta;k,\ell)}_{q,E}=\lim
_{n\to\infty}c'^{(\beta;k,\ell)}_{q,E}(n).\]
\end{defn}
The existence of the limits $c^{(\beta;k,\ell)}_{p,E}$ and
$c'^{(\beta;k,\ell)}_{q,E}$ is asserted by Theorem~\ref{ratio2}
below.

\medskip

 Definitions~\ref{def1} and~\ref{def??} obviously imply
\begin{rem}\label{hook1}

\[c^{(\beta;k,\ell)}_{p,E}(n)=\left ( \frac{\sum_{\lm\in H(k,\ell;n)}\lm_{p,n} \left (f^\lm\right)^\beta}
{\sum_{\lm\in H(k,\ell;n)}\left (
f^\lm\right)^\beta}-\frac{n}{k+\ell}\right )\cdot \frac{1}{\sqrt
n},\]

and
\[c'^{(\beta;k,\ell)}_{q,E}(n)=\left ( \frac{\sum_{\lm\in H(k,\ell;n)}\lm'_{q,n} \left (f^\lm\right)^\beta}
{\sum_{\lm\in H(k,\ell;n)}\left (
f^\lm\right)^\beta}-\frac{n}{k+\ell}\right )\cdot \frac{1}{\sqrt
n}.\]
\end{rem}
\begin{thm}\label{ratio2}
 Let $1\le p\le k$, $1\le q\le
\ell$, then the limits $c^{(\beta;k,\ell)}_{p,E}$ and
$c'^{(\beta;k,\ell)}_{q,E}$ exist, and are given as follows.
\begin{eqnarray}\label{hookeqn1}
c^{(\beta;k,\ell)}_{p,E}=\frac{\int_{\Omega_{k}} x_p\cdot \left
[D_k(x_1,\ldots,x_k) \cdot e^{-\frac{k+\ell}{2}(x_1^2+\cdots
+x_k^2)}\right ]^\beta d^{(k-1)}(x)}{\int_{\Omega_{k}} \left
[D_k(x_1,\ldots,x_k) \cdot e^{-\frac{k+\ell}{2}(x_1^2+\cdots
+x_k^2)}\right ]^\beta d^{(k-1)}(x)},
\end{eqnarray}
 and
\begin{eqnarray}\label{hookeqn2} c'^{(\beta;k,\ell)}
_{q,E}=\frac{\int_{\Omega_{\ell}} y_q\cdot \left [D_\ell
(y_1,\ldots,y_\ell) \cdot e^{-\frac{k+\ell}{2}(y_1^2+\cdots
+y_\ell ^2)}\right ]^\beta d^{(\ell -1)}(y)}{\int_{\Omega_{\ell}}
\left [D_\ell (y_1,\ldots,y_\ell) \cdot
e^{-\frac{k+\ell}{2}(y_1^2+\cdots +y_\ell ^2)}\right ]^\beta
d^{(\ell -1)}(y)}.\end{eqnarray} Thus, as $n$ goes to infinity,
\[\lm^{(\beta;k,\ell)}_{p,E}(n)\sim\frac{n}{k+\ell}+c^{(\beta;k,\ell)}_{p,E}\cdot\sqrt
n\quad\mbox{and}\quad
\lm'^{(\beta;k,\ell)}_{q,E}(n)\sim\frac{n}{k+\ell}+c'^{(\beta;k,\ell)}_{q,E}\cdot\sqrt
n.\]
\end{thm}

{\bf Proof}. We prove for $c^{(\beta;k,\ell)}_{p,E}$ -- in two
steps.

\smallskip

{\bf Step 1}.
We claim that the limit $c^{(\beta;k,\ell)}_{p,E}$ exists, and is
given by \begin{eqnarray}\label{eqnhook1}
c^{(\beta;k,\ell)}_{p,E}=\frac{\int_{\Omega_{(k,\ell)}} x_p\cdot
\left [D_k(x)\cdot D_\ell (y) \cdot e^{-\frac{k+\ell}{2}(\sum
x_i^2+\sum y_j^2)}\right ]^\beta d^{(k+\ell
-1)}(x;y)}{\int_{\Omega_{(k,\ell)}} \left [D_k(x)\cdot D_\ell (y)
\cdot e^{-\frac{k+\ell}{2}(\sum x_i^2+\sum y_j^2)}\right ]^\beta
d^{(k+\ell-1)}(x;y)}.
\end{eqnarray}
The proof of Equation~(\ref{eqnhook1}) is essentially the same as
that of Equation~(\ref{stripeqn1}). Its starting point is
Remark~\ref{hook1} (instead of Equation~(\ref{propdef1})). Here we
write, for $1\le i\le k$ and for $1\le j\le\ell$,
\[\lm_{i,n}=\frac{n}{k+\ell}+c_{i,n}\cdot\sqrt n\quad\mbox{ and}\quad
\lm'_{j,n}=\frac{n}{k+\ell}+c'_{j,n}\cdot\sqrt n,\] and consider
$\lm$--s with $c_i$ and $c'_j$ bounded in some interval. Now
follow a 'hook'--generalization of the proof of
Theorem~\ref{ratio}, applying Theorems~\ref{sheli3}
and~\ref{sheli4}, and complete the proof of
Equation~(\ref{eqnhook1}).

\medskip

{\bf Step 2}. We now transform (\ref{eqnhook1}) into
(\ref{hookeqn1}). Let $I^{(\beta;k,\ell)}_p$ denote the numerator
of (\ref{eqnhook1}):
\[I^{(\beta;k,\ell)}_p=\int_{ \Omega_{(k,\ell)} }x_p\left [D_k(x)\cdot D_\ell (y)\cdot e^{-((k+\ell)/2)
(\sum x_i^2+\sum y_j^2)}\right ]^\beta d^{(k+\ell-1)}(x;y)\,.\]

\medskip

Setting $\sum x_i=u$ we have $\sum y_j=-u$, and
$I^{(\beta;k,\ell)}_p=\int_{-\infty}^\infty K(u)L(-u)\,du$, where
\[K(u)= \int_{M_k(x,u)}x_p\left [ D_k(x)\cdot
e^{-((k+\ell)/2)(\sum x_i^2)}\right]^\beta d^{(k-1)}x\] and
\[L(-u)= \int_{M_{\ell}(y,-u)}\left [ D_\ell (y)\cdot
e^{-((k+\ell)/2)(\sum y_j^2)}\right]^\beta d^{(k-1)}x.\] Here
$M_k(x,u)=\{(x_1,\ldots,x_k)\mid x_1\ge \cdots\ge x_k
~\mbox{and}~\sum x_i=u\}$ and similarly,
$M_\ell(y,-u)=\{(y_1,\ldots,y_\ell)\mid y_1\ge\cdots\ge y_\ell~~\mbox{and}~ \sum y_j=-u\}$.\\

\smallskip

To evaluate $I^{(\beta;k,\ell)}_p$, proceed as follows. In $K(u)$
and $L(-u)$ substitute $x_i'=x_i-(u/k)$, and $y_j'=y_j+(u/\ell)$.
The Jacobians are $=1$,
~$x_t=x_t'+(u/k);~~D_k(x')=D_k(x);~~D_\ell(y')=D_\ell (y);~~\sum
x_i^2=\sum x_i'^2+(u^2/k)$ and $\sum y_j^2=\sum
y_j'^2+(u^2/\ell)$. Replacing $x_i'$ by $x_i$ and $y_j'$ by $y_j$,
it follows that
\[ I_p^{(\beta;k,\ell)}=J_1(x)\cdot J_3(y)\cdot A(u)
+ J_2(x)\cdot J_3(y)\cdot B(u).\]
Here \[A(u)=\int_{-\infty}^\infty
e^{-\frac{(k+\ell)^2u^2\beta}{2k\ell}} du,\qquad
B(u)=\int_{-\infty}^\infty \frac{u}{k}\cdot
e^{-\frac{(k+\ell)^2u^2\beta}{2k\ell}} du,\]
\[ J_1(x)=\int _{M_k(x,0)} x_p\cdot \left [D_k(x)\cdot
e^{-\frac{k+\ell}{2}(\sum x_i^2)}\right ]^\beta d^{(k-1)}x,\qquad
(J_1(x)=1\quad\mbox{if}\quad k=1),\]
\[ J_2(x)=\int _{M_k(x,0)}\left [D_k(x)\cdot
e^{-\frac{k+\ell}{2}(\sum x_i^2)}\right ]^\beta d^{(k-1)}x,\] and
\[J_3(y)=\int _{M_\ell (y,0)}\left [D_\ell(x)\cdot
e^{-\frac{k+\ell}{2}(\sum y_j^2)}\right ]^\beta d^{(\ell-1)}x,
\qquad (J_3(y)=1\quad\mbox{if}\quad \ell=1),
\] $(M_r(x,0)=\{(x_1,\ldots,x_r)\mid x_1\ge\cdots\ge x_r
~\mbox{and}~x_1+\cdots +x_r=0\}$, etc.). Since, trivially,
$B(u)=0$, deduce that  $I_p^{(\beta;k,\ell)}=J_1(x)\cdot
J_3(y)\cdot A(u)$.

\smallskip

Let  $\bar I^{(\beta;k,\ell)}$ denote the denominator in
Equation~(\ref{eqnhook1}). By exactly the same arguments it
follows that $\bar I^{(\beta;k,\ell)}=J_2(x)\cdot J_3(y)\cdot
A(u)$. By (\ref{eqnhook1})
\[c^{(\beta;k,\ell)}_{p,E}=\frac{I_p^{(\beta;k,\ell)}}{\bar
I^{(\beta;k,\ell)}}=\frac{J_1(x)}{J_2(x)},\] which is the
right-hand-side of Equation~(\ref{hookeqn1}). This completes the
proof.~~~~~~~~q.e.d.
\begin{rem}
The denominator-integral in~(\ref{eqnhook1}) is

\[\bar I^{(\beta;k,\ell)}={\int_{\Omega_{(k,\ell)}} \left
[D_k(x_1,\ldots,x_k)\cdot D_\ell (y_1,\ldots,y_\ell) \cdot
e^{-\frac{k+\ell}{2}(x_1^2+\cdots +x_k^2+y_1^2+\cdots
+y_\ell^2)}\right ]^\beta d^{(k+\ell-1)}(x;y)},\] and is
calculated explicitly in [~\cite{Be.Re.}, sec 7] (here
$\beta=2z$):
\[I^*(k,\ell)=I(k,\ell,2)=\frac{1}{k!\ell !}
\sqrt {2 \pi^{k+\ell -1}}\cdot\sqrt{\frac{\beta}{2\pi}} \cdot
\left (\frac {1}{\beta (k+\ell)}\right )^{\frac{1}{2}[(k(k-1)+\ell
(\ell-1))(\beta /2)+k+\ell]}~~~~~\]\[~~~~~~~~~~~~~~~~~~~~~\times
\frac{\prod_{i=1}^k\Gamma(i\beta/2+1)\cdot \prod_{j=1}^\ell
\Gamma(j\beta/2+1)}{\Gamma (\beta /2 +1)},\] where $\Gamma$ is the
{\it Gamma} function $(\Gamma (n+1)=n!)$.
\end{rem}
Theorems~\ref{ratio} and~\ref{ratio2} imply
\begin{cor}\label{relate1}
Let $1\le p\le k$ and $1\le q\le \ell$. Then
\[c^{(\beta;k,\ell)}_{p,E}=\sqrt{\frac{k}{k+\ell}}\cdot
c^{(\beta;k,0)}_{p,E},\qquad\mbox{and similarly} \qquad
c'^{(\beta;k,\ell)}_{q,E}=\sqrt{\frac{\ell}{k+\ell}}\cdot
c'^{(\beta;0,\ell)}_{q,E}.\]
\end{cor}
{\bf Proof.} Let $\alpha=\frac{k}{k+\ell}$ and in Equation
(\ref{hookeqn1}) substitute $x=\sqrt\alpha\cdot v$. By routine
calculations, this substitution transforms the ration of integrals
(\ref{hookeqn1}) into the ratio in (\ref{stripeqn1}) -- multiplied
by the factor $\sqrt\alpha$, which completes the proof.

\medskip

\section{The distribution functions }\label{sec6}
\subsection{First term approximation}
\begin{defn}\label{def111}
Let $z>0$. Denote
\[H(k,\ell;n,z)=\{\lm\in H(k,\ell;n)\mid
 \lm_{1,n},\;\lm'_{1,n}\le\frac{n}{k+\ell}+z\sqrt n\,\}.\]

Let $\beta >0$, $1\le p\le k$ and $1\le q\le\ell$. The
distribution of the length of the $p$--th row as a function of $z$
is defined  as
\begin{eqnarray}\label{od1}
\lm^{(\beta;k,\ell)} _p(n,z)=\frac{\sum_{\lm\in H(k,\ell;n,z)
 }\lm_{p,n}
\cdot(f^\lm)^{\beta}} {\sum_{\lm\in H(k,\ell;n)
}(f^\lm)^{\beta}}.\end{eqnarray}

Similarly, the distribution of the length of the $q$--th column as
a function of $z$ is defined  as
\begin{eqnarray}\label{od123}
\lm'^{(\beta;k,\ell)} _q(n,z)=\frac{\sum_{\lm\in H(k,\ell;n,z)
 }\lm'_{q,n}
\cdot(f^\lm)^{\beta}} {\sum_{\lm\in H(k,\ell;n)
}(f^\lm)^{\beta}}.\end{eqnarray}

\end{defn}
Let $\lm\in H(k,\ell;n,z)$ with $n$ large, then necessarily
$\lm_1,\ldots,\lm_k\ge \ell$ and $\lm'_1,\ldots,\lm'_{\ell}\ge k$
(in fact, Lemma~\ref{ch} proves a much stronger property) so that
$\lm_1,\ldots,\lm_k+\lm'_1,\ldots,\lm'_{\ell}=n+k\ell$. Write
\[\lm_p=\frac{n+k\ell}{k+l}+c_p\cdot\sqrt n,~~p=1,\ldots,k\quad\mbox{
and}\quad \lm'_q=\frac{n+k\ell}{k+l}+c'_q\cdot\sqrt n,~~
q=1,\ldots,\ell,\] and notice that $\sum c_p+\sum c'_q=0$.

\begin{lem}\label{ch}
With the above notations (and $n$ large),
\[-\frac{\ell+p-1}{k-p+1}\cdot z\le c_p\le z,~~p=1,\ldots,k
\qquad\mbox{and}\qquad -\frac{k+q-1}{\ell-q+1}\cdot z\le c'_q\le
z,~~q=1,\ldots,\ell .\]
\end{lem}
{\bf Proof.} Clearly, all $c_p,\;c'_q\le z$. Assume for example
that \[ c_p  < -\frac{\ell+p-1}{k-p+1}\cdot z,\quad\mbox{hence
also}\quad c_k,c_{k-1},\ldots,c_p < -\frac{\ell+p-1}{k-p+1}\cdot
z.\qquad \mbox{Thus}\]
\[0=c_1+\cdots +c_k+c'_1+\cdots
+c'_{\ell}<-(k-p+1)\cdot\frac{\ell+p-1}{k-p+1}\cdot z+c_1+\cdots
+c_{p-1}+c'_1+\cdots +c'_{\ell}~~~~~~~~~\]
\[~~~~~~~~~~~~~~~~~~~~~~~~~~~~~~~~~~~~~~~~~~~~\le -(\ell+p-1)\cdot
z+\ell+p-1)\cdot z=0,\]  a contradiction. Similarly for $c'_q$.
This proves the lemma.

\bigskip

Recall that\\
$\Omega_{(k,\ell)}=\{(x_1,\ldots,x_k,y_1,\ldots,y_\ell)\mid
x_1\ge\cdots\ge x_k\,;~~ y_1\ge\cdots\ge y_\ell\,;~~\sum x_i+\sum
y_j=0\}$ and
$\Omega_{(k,\ell),z}=\{(x_1,\ldots,x_k;y_1,\ldots,y_\ell)\in\Omega_{(k,\ell)}\mid
x_1,\;y_1\le z\}$.
For example, $\Omega_{(2,0),z}=\{(x,-x)\mid 0\le x\le z\}$.

\begin{thm}\label{first222} Let $\beta,\,z > 0$ and
denote
\begin{eqnarray}\label{dist222} r_{(k,\ell),\beta}(z)=
\frac{\int_{\Omega_{(k,\ell),z}} \left [ D_k(x)\cdot
D_\ell(y)\cdot e^{-\frac{k+\ell}{2}(\sum x_i^2+\sum y_j^2)}\right
]^\beta d^{(k+\ell-1)}(x;y)}{\int_{\Omega_{(k,\ell)}} \left
[D_k(x)\cdot  D_\ell(y) \cdot e^{-\frac{k+\ell}{2}(\sum x_i^2+\sum
y_j^2)}\right ]^\beta d^{(k+\ell-1)}(x;y)}.\end{eqnarray}  Let
$1\le p\le k$, $1\le q\le\ell$ and let $n\to\infty$, then
\[\lm^{(\beta;k,\ell)}_{p}(n,z),\;\lm'^{(\beta;k,\ell)}_{q}(n,z)
\sim\frac{n}{k+\ell}\cdot r_{(k,\ell),\beta}(z).\]

In other words,
\[\lim_{n\to\infty}\frac{k+\ell}{n}\cdot\lm^{(\beta;k,\ell)}
_p(n,z) =\lim_{n\to\infty}\frac{k+\ell}{n}\cdot\frac{\sum_{\lm\in
H(k,\ell;n,z)
 }\lm_{p,n}
\cdot(f^\lm)^{\beta}} {\sum_{\lm\in H(k,\ell;n)
}(f^\lm)^{\beta}}=\]
\begin{eqnarray}\label{od22}
=\frac{\int_{\Omega_{(k,\ell),z}} \left [ D_k(x)\cdot D_\ell(y)
\cdot e^{-\frac{k+\ell}{2}(\sum x_i^2+\sum y_j^2)}\right ]^\beta
d^{(k+\ell-1)}(x;y)}{\int_{\Omega_{(k,\ell)}} \left [D_k(x)\cdot
D_\ell(y) \cdot e^{-\frac{k+\ell}{2}(\sum x_i^2+\sum y_j^2)}\right
]^\beta d^{(k+\ell-1)}(x;y )}.\end{eqnarray} Similarly for
\[\lim_{n\to\infty}\frac{k+\ell}{n}\cdot\lm'^{(\beta;k,\ell)}
_q(n,z).\]
\end{thm}
{\bf Proof} is similar to the proof of Theorem~\ref{first}: Let
$\lm_{p,n}=\frac{n}{k+\ell}+c_{p,n}\cdot\sqrt n$, then
\[\frac{k+\ell}{n}\cdot\lm^{(\beta;k,\ell)}
_p(n,z)= \frac{\sum_{\lm\in H(k,\ell;n,z)
 }
(f^\lm)^{\beta}} {\sum_{\lm\in H(k,\ell;n)
}(f^\lm)^{\beta}}+\frac{k+\ell}{\sqrt n}\cdot \left
(\frac{\sum_{\lm\in H(k,\ell;n,z)
 }c_{p,n}
\cdot(f^\lm)^{\beta}} {\sum_{\lm\in H(k,\ell;n) }(f^\lm)^{\beta}}
\right ).\]

By Theorem~\ref{sheli5}, the first summand approaches $
r_{(k,\ell),\beta}(z)$ as $n$ goes to infinity, and by
Lemma~\ref{ch}, $|c_{p,n}|\le b\cdot z$ for an appropriate
constant $b>0$, therefore the second summand obviously goes to
zero as $n$ goes to infinity.~~~~~~~~~~q.e.d.

\medskip

Theorem \ref{first222} is a first--term approximation of the
distribution function.

\subsection{Conjectures about the second term
approximation}\label{sec7}

Let \begin{eqnarray}\label{od11}
s_{(k,\ell),\beta}(n,z)=\frac{\sum_{\lm\in H(k,\ell;n,z)
 }(f^\lm)^{\beta}} {\sum_{\lm\in H(k,\ell;n)
}(f^\lm)^{\beta}}.\end{eqnarray} As in the proof of
Theorem~\ref{first222}, $~\lim_{n\to\infty}s_{(k\ell),\beta}(n,z)=
r_{(k,\ell),\beta}(z)$. Numerical evidence suggest the following
(vague) conjecture.
\begin{conj}\label{con1}
For all $k,\,\ell\ge 0$ and $\beta,\, z>0$, as $n$ goes to
infinity the expression
\begin{eqnarray}\label{od2}
\frac{\sqrt n}{k+\ell}\cdot\left[
s_{(k,\ell),\beta}(n,z)-r_{(k,\ell),\beta}(z)\right
]=~~~~~~~~~~~~~~~~~~~~~~~~~~~~~~~~~~~~~~~~~~~~~~~~~~~~~~~
\end{eqnarray}
\[= \frac{\sqrt n}{k+\ell}\cdot\left[ \frac{\sum_{\lm\in
H(k,\ell;n,z)
 }(f^\lm)^{\beta}} {\sum_{\lm\in H(k,\ell;n)
}(f^\lm)^{\beta}} -\frac{\int_{\Omega_{(k,\ell),z}} \left [
D_k(x)\cdot D_\ell(y) \cdot e^{-\frac{k+\ell}{2}(\sum x_i^2+\sum
y_j^2)}\right ]^\beta
d^{(k+\ell-1)}(x;y)}{\int_{\Omega_{(k,\ell)}} \left [D_k(x)\cdot
D_\ell(y) \cdot e^{-\frac{k+\ell}{2}(\sum x_i^2+\sum y_j^2)}\right
]^\beta d^{(k+\ell-1)}(x )} \right ]\] oscillates in some
symmetric bounded interval centered at zero. We denote that
interval as $(-L(k,\ell,\beta,z),L(k,\ell,\beta,z))$, where (we
conjecture that) $-L(k,\ell,\beta,z)$ and $L(k,\ell,\beta,z)$ are
the respective infimum and supremum of the values in Equation
(\ref{od2}).
\end{conj}

\begin{defn}\label{def???2}
Let $\beta,\,z>0$, $1\le p\le k$ and $1\le q\le\ell$. Then
$\lm^{(\beta;k,\ell)}_p(n,z)$ is given by Equation~(\ref{od1}).
Similarly for the columns $\lm'^{(\beta;k,\ell)}_q(n,z)$. Now
define $c^{(\beta;k,\ell)}_p(n,z)$ via the equation
\[\lm^{(\beta;k,\ell)}_{p}(n,z)=\frac{n}{k+\ell}\cdot
r_{(k,\ell),\beta}(z)+c^{(\beta;k,\ell)}_{p}(n,z)\sqrt n\,.\]
\end{defn}
Similarly for  $c'^{(\beta;k,\ell)}_q(n,z)$. We would like to
understand the behavior of $c^{(\beta;k,\ell)}_{p}(n,z)$ and
$c'^{(\beta;k,\ell)}_q(n,z)$ as $n$ goes to infinity. We consider
$c^{(\beta;k,\ell)}_{p}(n,z)$.
\medskip

 Note that Definition~\ref{def111}, Theorem \ref{first222} and Definition~\ref{def???2} imply

\begin{prop}\label{hook111}
\begin{eqnarray}\label{eq4}
c^{(\beta;k,\ell)}_{p}(n,z)=\left ( \frac{\sum_{\lm\in
H(k,\ell;n,z)}\lm_{p,n} \left (f^\lm\right)^\beta} {\sum_{\lm\in
H(k,\ell;n)}\left ( f^\lm\right)^\beta}-\frac{n}{k+\ell}\cdot
r_{(k,\ell),\beta}(z)\right )\cdot \frac{1}{\sqrt
n}.\end{eqnarray}\end{prop}
\begin{conj}\label{con2}
Recall the interval $(-L(k,\ell,\beta,z),L(k,\ell,\beta,z))$ from
Conjecture~\ref{con1}. As $n$ goes to infinity,
$c^{(\beta;k,\ell)}_{p}(n,z)$ oscillates in the interval
$(-L(k,\ell,\beta,z),L(k,\ell,\beta,z))+ s_p^{(k,\ell),\beta}(z)$,
where
\begin{eqnarray}\label{eqnhook111} s_p^{(k,\ell),\beta}(z)=
\frac{\int_{\Omega_{(k,\ell),z}} x_p\cdot \left [D_k(x)\cdot
D_\ell(y)\cdot e^{-\frac{k+\ell}{2}(\sum x_i^2+\sum y_j^2)}\right
]^\beta d^{(k+\ell -1)}(x)}{\int_{\Omega_{(k,\ell)}} \left
[D_k(x)\cdot D_\ell(y)\cdot e^{-\frac{k+\ell}{2}(\sum x_i^2+\sum
y_j^2)}\right ]^\beta d^{(k+\ell-1)}(x)}.
\end{eqnarray}
\end{conj}
{\bf The proof} is based on Conjecture~\ref{con1} as follows. In
(\ref{eq4}) write $\lm_{p,n}=\frac{n}{k+\ell}+c_{p,n}\cdot\sqrt
n$, then $c^{(\beta;k,\ell)}_{p}(n,z)= A(n)+B(n)$, where
\[A(n)=\frac{\sqrt n}{k+\ell}\cdot\left ( \frac{\sum_{\lm\in H(k,\ell;n,z)}
\left (f^\lm\right)^\beta} {\sum_{\lm\in H(k,\ell;n)}\left (
f^\lm\right)^\beta}-r_{(k,\ell),\beta}(z)\right )\] and
\[B(n)=\frac{\sum_{\lm\in H(k,\ell;n,z)} c_{p,n}\cdot\left (f^\lm\right)^\beta}
{\sum_{\lm\in H(k,\ell;n)}\left ( f^\lm\right)^\beta}.\] By
arguments similar to those in previous proofs,
$\lim_{n\to\infty}B(n)=s_{(k,\ell),\beta}(z)$, and by
Conjecture~\ref{con1}, $A(n)$ oscillates in the interval
$(-L(k,\ell,\beta,z),L(k,\ell,\beta,z))$.

\part{Some special cases}
\section{$\beta=2$, comparison of expected and maximal
shapes}\label{sec8}
\subsection{Expected shape}
Recall the RSK bijection $\sg\longleftrightarrow (P_\lm,Q_\lm )$
and the subsets
\[S_{k,\ell;n}=\{\sg\in S_n\mid\mbox{under the RSK}~~\sg\longleftrightarrow
(P_\lm,Q_\lm),\quad \lm\in  H(k,\ell;n)\,\}\] from
Section~\ref{rsk}. For example if $\ell=0$,  it follows from well
known properties of the RSK that $S_{k,0;n}$ are the permutations
in $S_n$ with longest decreasing subsequence having length $\le
k$. In general, if $\sg \in S_{n}$ is of shape
$\lm=(\lm_1,\lm_2,\ldots)$, then $\lm_1$ is the length of a
maximal increasing subsequence in $\sg$. Thus, for example,
$\lm^{(2;k,\ell)}_{1,E}$ (see Definition~\ref{def1}) is the
expected length of the longest increasing subsequences in
$S_{k,\ell;n}$.

\medskip

Consider first the case $\ell=0$. Since the case $k=1$ is trivial,
assume $k\ge 2$ (and $\ell=0$). In that case, by
Theorem~\ref{ratio}, the expected shape in $S_{k,0;n}$ is\\
$\lm^{(2;k,0)}_{E}\sim \left (\frac{n}{k}+c_1\cdot\sqrt
n\,,\,\frac{n}{k}+c_2\cdot\sqrt
n\,,\ldots,\frac{n}{k}+c_k\cdot\sqrt n\right )$, where for $1\le
p\le k$,
\begin{eqnarray}\label{beta2}
~~c_{p}=c^{(2;k,0)}_{p,E}=\frac{\int_{\Omega_k} x_p\cdot \left
[D_k(x_1,\ldots,x_k)\cdot e^{-\frac{k}{2}(x_1^2+\cdots
+x_k^2)}\right ]^2 d^{(k-1)}x}{\int_{\Omega_k} \left
[D_k(x_1,\ldots,x_k)\cdot e^{-\frac{k}{2}(x_1^2+\cdots
+x_k^2)}\right ]^2 d^{(k-1)}x}. \end{eqnarray} In general, the
expected shape of the permutations in $S_{k,\ell;n}$ is given by
Definition~\ref{def1}, and the asymptotic shape as $n\to\infty$ is
given by Theorem~\ref{ratio2}, both with $\beta=2$.

\subsection{Maximal $f^\lm$}\label{max5}

Given a subset of partitions $\Gamma_n\subseteq Y_n$, one looks
for $\lm=\lm_{max}\in\Gamma_n$ with maximal degree $f^\lm$, see
Section~\ref{1.2}.
In the case $\Gamma_n=H(k,0;n)$, the asymptotics of $\lm_{max}$ is
given in \cite{AR}, which we briefly describe here. The analogue
result from \cite{hook}, for the case $\Gamma_n=H(k,\ell;n)$, is
also described below.

\smallskip

Let $H_k(x)$ denote the $k$-th Hermit polynomial. It is defined
via the equation \[\frac{d^k}{dx^k}\left (e^{-x^2} \right
)=(-1^k)H_k(x)e^{-x^2}. \] Thus
$H_0(x)=1,~H_1(x)=2x,~H_2(x)=4x^2-2$, ~$H_3(x)=4x(2x^2-3)$,
~$H_4(x)=16x^4-48x^2+12$, etc. $H_k(x)$ is of degree $k$ and its
roots are real and distinct, denoted
\[x_1^{(k)}< x_2^{(k)}<\cdots< x_k^{(k)}.\]
Also,  $~x_1^{(k)}+ x_2^{(k)}+\cdots+ x_k^{(k)}=0$ . The following
theorem is proved in \cite{AR}:

\begin{thm}
As $n\to\infty$, the maximum $~\max\{f^\lm\mid\lm\in H(k,0;n)\}$
occurs when
\[\lm\sim\lm^{(k,0)}_{max}
=\left ( \frac{n}{k}+x_k^{(k)}\sqrt {\frac{k}{
n}}\,,\ldots,\frac{n}{k}+x_1^{(k)}\sqrt {\frac{k}{ n}} \right ).\]
\end{thm}

The analogue result for $\Gamma_n=H(k,\ell ;n)$ is given in
\cite{hook}:
\begin{thm}
Let $\lm=\lm_{max}^{(k,\ell)}\in H(k,\ell;n)$ maximize $f^\lm$ in
$H(k,\ell;n)$ and write
$\lm_{max}^{(k,\ell)}=(\lm_1,\ldots,\lm_k,\ldots)$ ~and
$\lm_{max}^{'(k,\ell)}=(\lm'_1,\ldots,\lm'_\ell,\ldots)$, ~and
assume $n\to\infty$, then
\[(\lm_1,\ldots,\lm_k)\sim\left(\frac{n}{k+\ell}+x^{(k)}_k\sqrt\frac{n}{k+\ell}
\,\,,\ldots,\frac{n}{k+\ell}+x^{(k)}_1\sqrt\frac{n}{k+\ell}\;\right)\]
and
\[(\lm'_1,\ldots,\lm'_\ell)\sim\left(\frac{n}{k+\ell}+x^{(\ell)}_\ell\sqrt\frac{n}{k+\ell}
\,\,,\ldots,\frac{n}{k+\ell}+x^{(\ell)}_1\sqrt\frac{n}{k+\ell}\;\right).\]
Here $x^{(k)}_1<\cdots <x^{(k)}_k$ are the roots of the $k$--th
Hermit polynomial and similarly for the $x^{(\ell)}_j$--s.
\end{thm}

\subsection{Examples of some $(k,0)$ cases}\label{ex4}

Let $\Gamma_n=H(k,0;n)$ and $\beta=2$~($\longleftrightarrow_{RSK}
S_{k,0;n}$) and compare the expected shape $\lm_E^{(2;k,0)}$ with
the maximizing shape $\lm_{max}^{(k,0)}$. We begin with

\medskip


\medskip

{\bf The case $k=2$}.\\ Here the expected shape is given by
$\lm^{(2;2,0)}_{E}\sim\left (\frac{n}{2}+c^{(2;2,0)}_{1,E}\sqrt
n\,,\frac{n}{2}+c^{(2;2,0)}_{2,E}\sqrt n\right )$, where
$c^{(2;2,0)}_{2,E}=-c^{(2;2,0)}_{1,E}$, and by
Equation~(\ref{beta2}), ~~$c^{(2;2,0)}_{1,E}=I_1/I_2$ is the ratio
of the following integrals:
\[I_1=\int_0^{\infty} x\left [
2x\cdot e^{-2x^2}\right ]^2\,dx=\frac{1}{8},\qquad\mbox{
and}\qquad I_2=\int_0^{\infty}\left [ 2x\cdot e^{-2x^2}\right
]^2\,dx=
\frac{\sqrt\pi}{8}.\] Thus $c^{(2;2,0)}_{1,E}=\frac{1}{\sqrt\pi},$
hence
\begin{eqnarray}\label{case2}\lm^{(2;2,0)}_{E}\sim \left
(\frac{n}{2}+\frac{1}{\sqrt\pi}\cdot\sqrt
n\;,\;\frac{n}{2}-\frac{1}{\sqrt\pi}\cdot\sqrt n\right )= \left
(\frac{n}{2}+0.56419\sqrt n\;,\;\frac{n}{2}-0.56419\sqrt n\right
).\end{eqnarray}

Compare $\lm^{(2;2,0)}_{E}$ with $\lm^{(2,0)}_{max}$: Here
$x_2^{(2,0)}=1/\sqrt 2,~x_1^{(2,0)}=-1/\sqrt 2$, so
\[\lm^{(2,0)}_{max}=\left(\frac{n}{2}+\frac{1}{2}\sqrt
n\;,\;\frac{n}{2}-\frac{1}{2}\sqrt
n\right)=\left(\frac{n}{2}+0.5\sqrt n\;,\;\frac{n}{2}-0.5\sqrt
n\right).\]
{\bf Note:} working with Equation~(\ref{propdef1}) and with
'Mathematica' we calculated $c^{(2;2,0)}_{1,E}(n)$. For
$n=100,~200,~ 300,~400,~500,~600,~700,~800,~900$, the
corresponding values of $c^{(2;2,0)}_{1,E}(n)$ are:\\ $0.517699,
~0.530593, ~0.536496, ~0.54007, ~0.542533, ~0.544364, ~0.545795,
~0.546952$ and $0.547915$, agreeing with Theorem~\ref{ratio}.

\bigskip

{\bf The case $k=3$}

\medskip

Here $\lm^{(2;3,0)}_{E}\sim \left(
\frac{n}{3}+c^{(2;3,0)}_{1,E}\sqrt
n\;,\;\frac{n}{3}+c^{(2;3,0)}_{2,E}\sqrt
n\;,\;\frac{n}{3}+c^{(2;3,0)}_{3,E}\sqrt n \right )$, so we
calculate the $c$--s. By Equation~(\ref{beta2}),
$c^{(2;3,0)}_{1,E}=J_1/J_2$ is the ratio of the following
integrals:
\[J_1=\int_{\Omega_3}x_1\left[D_3(x)e^{-\frac{3}{2}(x_1^2+x_2^2+x_3^2)}\right]^2\,
d^{(2)}x~\mbox{and}~ J_2=\int_{\Omega_3}\left
[D_3(x)e^{-\frac{3}{2}(x_1^2+x_2^2+x_3^2)}\right]^2\,d^{(2)}x.
\] By
Remark~\ref{value1} with $\beta=2$ and $k=3$,
$J_2=\frac{\pi}{324\sqrt 3}$.
 We
calculate the numerator $J_1$.
The domain $\Omega_3$ of integration is defined by: $x_1\ge x_2\ge
x_3$ and $x_3=-(x_1+x_2)$, so $x_1\ge x_2\ge -(x_1+x_2).$ When
$x_2\le 0$, that condition is equivalent to
$x_2\ge -x_1/2$. When $x_2\ge 0$, that condition is equivalent to
$x_2\le x_1$. It follows that
\[J_1=\int_0^{\infty} \left [\int_{\frac{-x_1}{2}}^{x_1}\left (
x_1\left
[(x_1-x_2)(2x_1+x_2)(x_1+2x_2)e^{-3(x_1^2+x_2^2+x_1x_2)}\right
]^2\right ) dx_2\right ]dx_1. \]

After some routine calculations ('Mathematica' was used here) we
obtain $J_1=\frac{\sqrt\pi}{288\sqrt 2}$.

It follows that when $k=3$,
\begin{eqnarray}\label{limit1} \lim_{n\to\infty}c^{(2;3,0)}_{1,E}(n)=
c^{(2;3,0)}_{1,E}=\frac{9\sqrt 3}{8\sqrt
{2\pi}}=0.777362...\end{eqnarray} By similar calculations it
follows that $c^{(2;3,0)}_{2,E}=0$, hence $
c^{(2;3,0)}_{3,E}=-c^{(2;3,0)}_{1,E}$. Thus
\[\lm^{(2;3,0)}_{E} \sim \left (\frac{n}{3}+\frac{9\sqrt 3}{8\sqrt
{2\pi}}\sqrt n\;,\;\frac{n}{3}\;,\;\frac{n}{3}-\frac{9\sqrt
3}{8\sqrt {2\pi}} \sqrt n \right )
=~~~~~~~~~~~~~~~~~~~~~~~~~~~~~~~~~~\]
\[~~~~~~~~~~~~~~~~~~~~~~~~~~~~~~~~~~~=\left (\frac{n}{3}+0.777362\sqrt
n\;,\;\frac{n}{3}\;,\;\frac{n}{3}-0.777362 \sqrt n \right ).\]

\medskip

Compare now with $\lm^{(3,0)}_{max}$. Since $H_3(x)=4x(2x^2-3),
~x_3^{(3,0)}=\frac{\sqrt 3}{\sqrt 2},~x_2^{(3,0)}=0$ and
$x_1^{(3,0)}=-\frac{\sqrt 3}{\sqrt 2}$. Thus
\[\lm^{(3,0)}_{max}\sim \left (\frac{n}{3}+\frac{1}{\sqrt
2}\sqrt n\; ,\;\frac{n}{3}\;,\;\frac{n}{3}-  \frac{1}{\sqrt
2}\sqrt n  \right )=\left (\frac{n}{3}+0.707107\sqrt n
\;,\;\frac{n}{3}\;,\;\frac{n}{3}- 0.707107\sqrt n  \right ).\]

{\bf Note:} For $n=200,~300,~400,~500,~700,~850$ and $1100$,
'Mathematica' and Equation~(\ref{propdef1}) give the following
corresponding values of $c^{(2;3,0)}_{1,E}(n)$: $0.719084,\;
0.729014,\; 0.735096,\; 0.739317,\; 0.74494,\; 0.747816$ and
$0.751261$, in accordance with Theorem~\ref{ratio}.

\section{Examples of some $(k,\ell)$--hook cases, $\beta=2$}\label{sec9}

\subsection{The $(1,1)$ case}
By Equation~(\ref{hookeqn1}) $c^{(2;1,1)}_{1,E}=0$.
[Alternatively, calculate $c^{(\beta;1,1)}_{1,E}$ by
Equation~(\ref{eqnhook1}). Here $\Omega_{1,1}=\{(x,y)\mid
x+y=0\}$, so $y=-x$ and $-\infty\le x\le\infty$. Thus
$c^{(\beta;1,1)}_{1,E}=I_1(\beta)/I_2(\beta)$ where
$I_1(\beta)=\int_{-\infty}^\infty x\left[
e^{-2x^2}\right]^\beta\,dx=0$.] Similarly $c'^{(2;1,1)}_{1,E}=0$
by Equation~(\ref{hookeqn2}). It follows that (for any $\beta >0$)
\[\lm^{(\beta;1,1)}_{E}\sim \left(\frac{n}{2},1^{\frac{n}{2}}\right).\]

\bigskip


\subsection{The $(2,1)$ case}
By Corollary~\ref{relate1}
\begin{eqnarray}\label{limit2} \lim_{n\to\infty}c^{(2;2,1)}_{1,E}(n)=c^{(2;2,1)}_{1,E}=\sqrt{\frac{2}{3}}\cdot
c^{(2;2,0)}_{1,E}=\sqrt{\frac{2}{3}}\cdot\frac{1}{\sqrt\pi}=0.460659
.\end{eqnarray} Since $c^{(2;2,0)}_{2,E}=-c^{(2;2,0)}_{1,E}$,
deduce that
\begin{eqnarray}\label{limit3} c^{(2;2,1)}_{2,E}=
-\sqrt{\frac{2}{3}}\cdot\frac{1}{\sqrt\pi}=-0.460659.\end{eqnarray}
Since the sum of the coordinates in $\lm^{(2;2,1)}_{E}(n)$ is $n$,
it follows that $c'^{(2;2,1)}_{1,E}=0$ (this can also be deduced
directly from Theorem \ref{ratio2}). Therefore
\[\lm^{(2;2,1)}_{E}\sim \left (\frac{n}{3}+0.460659\sqrt n\,,
~~\frac{n}{3}-0.460659\sqrt n\,,~~1^{n/3}\right ).\]

\medskip

{\bf Note:}  for $n=100,\;200,\;300,\;400,\;700$ and
$1100$,~Remark~(\ref{hook1}) and 'Mathematica' give the following
values of $c^{(2;2,1)}_{1,E}(n)$: $0.45423,\;0.45571,\;0.456475,\;
0.456962,\; 0.457777$ and $0.458317$, agreeing with
Theorem~\ref{ratio2}.




\section{Expected shape, $\beta=1$}\label{sec10}
\subsection{The general case}
Since $\sg\in S_n$ is an involution if''f the RSK yields
$\sg\longleftrightarrow (P_\lm,P_\lm )$, therefore
$\lm^{(1;k,\ell)}_{E}(n)$ is  the {\it expected shape of the
involutions in $S_{k,\ell;n}$}. Note that when $\Gamma_n=Y_n$ (and
$\beta=1$) Baik and Rains \cite{Ba.Ra.} showed that the expected
length of the first row (or of the longest increasing subsequence
in involutions in $S_n$) is again $2\sqrt n$, i.e.~the same as the
first row of $\lm_{max}$, ~see \cite{Ba.Ra.} (1.5). As we show
below, this is {\it not} the case when $\Gamma_n=H(k,\ell;n)$.

\medskip

Denote $\tilde\lm^{(k,\ell)}(n)=\lm^{(1;k,\ell)}_{E}(n)$. We
summarize:

\begin{enumerate}
\item
Let $1\le p\le k,~~1\le q\le \ell$, then
\[\tilde \lm^{(k,\ell)} _{p,n}=\frac{\sum_{\lm\in H(k,\ell;n)}\lm_{p,n}\cdot f^\lm}
{\sum_{\lm\in H(k,\ell;n)}f^\lm}\qquad\mbox{and}\qquad\tilde
\lm'^{(k,\ell)} _{q,n}=\frac{\sum_{\lm\in
H(k,\ell;n)}\lm'_{q,n}\cdot f^\lm} {\sum_{\lm\in
H(k,\ell;n)}f^\lm} .\]
\item
 Define $\tilde c^{(k,\ell)}_{p,n}$ via
 \[\tilde\lm^{(k,\ell)}_{p,n}=\frac{n}{k+\ell}+{\tilde c^{(k,\ell)}_{j,n}}\sqrt
 n,\qquad\mbox{and}\qquad
 \tilde c^{(k,\ell)}_{p}=\lim _{n\to\infty}\tilde c^{(k,\ell)}_{p,n}.\]
 Similarly for $\tilde c'^{(k,\ell)}_{q,n}~$ and $~\tilde c'^{(k,\ell)}_{q}$.
Thus, when $n\to\infty$,
\[\tilde\lm^{(k,\ell)}_{p,n}\sim\frac{n}{k+\ell}+\tilde c^{(k,\ell)}_p\cdot \sqrt
n\qquad\mbox{and}\qquad\tilde\lm'^{(k,\ell)}_{q,n}\sim\frac{n}{k+\ell}+\tilde
c'^{(k,\ell)}_q\cdot \sqrt n.\]
\item
(Theorem~\ref{ratio2}, $\beta=1$)   Let $1\le p\le k$, $1\le q\le
\ell$, then the limits $\tilde c^{(k,\ell)}_{p}$ and $\tilde
c'^{(k,\ell)}_{q}$ exist, and are given as follows.
\begin{eqnarray}\label{hookeqn3}
\tilde c^{(k,\ell)}_{p}=\frac{\int_{\Omega_{k}} x_p\cdot
D_k(x_1,\ldots,x_k) \cdot e^{-\frac{k+\ell}{2}(x_1^2+\cdots
+x_k^2)} d^{(k-1)}(x)}{\int_{\Omega_{k}} D_k(x_1,\ldots,x_k) \cdot
e^{-\frac{k+\ell}{2}(x_1^2+\cdots +x_k^2)}d^{(k-1)}(x)}.
\end{eqnarray}
Similarly for the expected
columns:\begin{eqnarray}\label{hookeqn4} \tilde
c'^{(k,\ell)}_{q}=\frac{\int_{\Omega_{\ell}} y_q\cdot D_\ell
(y_1,\ldots,y_\ell) \cdot e^{-\frac{k+\ell}{2}(y_1^2+\cdots
+y_\ell ^2)}d^{(\ell -1)}(y)}{\int_{\Omega_{\ell}} D_\ell
(y_1,\ldots,y_\ell) \cdot e^{-\frac{k+\ell}{2}(y_1^2+\cdots
+y_\ell ^2)}d^{(\ell -1)}(y)}.\end{eqnarray}
\item
(Corollary~\ref{relate1}, $\beta=1$)   We have
\[\tilde c^{(k,\ell)}_{p}=\sqrt{\frac{k}{k+\ell}}\cdot
\tilde c^{(k,0)}_{p},\qquad\mbox{and similarly}\qquad  \tilde
c'^{(k,\ell)}_{q}=\sqrt{\frac{\ell}{k+\ell}}\cdot \tilde
c'^{(0,\ell)}_{q}.\]
\end{enumerate}

\subsection{Examples for involutions in $S_{k,0;n}$}
(i.e.~$\beta=1$ and $\ell=0$).\\

{\bf When $k=2$}, we obtain $\int_0^{\infty}2x\cdot
e^{-2x^2}=\frac{1}{2}\quad\mbox {and}\quad
\int_0^{\infty} x2x\cdot e^{-2x^2}=\frac{\sqrt\pi}{4\sqrt 2}$.\\
Thus $\tilde c^{(2,0)}_1=\frac{\sqrt\pi}{2\sqrt 2}$ and
\[\tilde\lm^{(2,0)} =\left ( \tilde\lm^{(2,0)}_1\;,\tilde\lm^{(2,0)}_2\right )\sim\left
(\frac{n}{2}+  \frac{\sqrt\pi}{2\sqrt 2}\cdot\sqrt
n\;,\;\frac{n}{2}-\frac{\sqrt\pi}{2\sqrt 2} \cdot\sqrt n\right
)=~~~~~~~~~~~~~~~~~~~~~~~~~~~~~~~~~~~~~~\]

\begin{eqnarray}\label{inv1} =\left (\frac{n}{2}+ 0.626657\sqrt
n\;,\;\frac{n}{2}-0.626657\sqrt n\right ).\end{eqnarray}

\medskip

{\bf Note:} for $n=200,\,400,\,600,\,800,\,1000,\,1200$ and
$1400$, Remark~\ref{hook1} and 'Mathematica' give the following
values of $\tilde c^{(2,0)}_{1}(n)$:
$0.592086,\,0.602049,\,0.606506$, ~$0.609175$, $0.611002$, $
0.612354 $ and $0.613406$, agreeing with
Equation~(\ref{hookeqn3}).

\bigskip

{\bf When $k=3$} we have

By Equation~(\ref{hookeqn3}), $c^{(1;3,0)}_{1,E}={\tilde
J_1}/{\tilde J_2}$ is the ratio of the following integrals:
\[\tilde J_1=\int_{\Omega_3}x_1D_3(x)e^{-\frac{3}{2}(x_1^2+x_2^2+x_3^2)}\,
d^{(2)}x\quad\mbox{and}\quad \tilde J_2=\int_{\Omega_3}
D_3(x)e^{-\frac{3}{2}(x_1^2+x_2^2+x_3^2)}\,d^{(2)}x
\] By
Remark~\ref{value1}  ($\beta=1$,  $k=3$), $\tilde J_2=\frac{\sqrt
\pi}{27}$. We calculate the numerator $\tilde J_1$. Similar to the
evaluation of $J_1$ in Section \ref{ex4},  here
\[\tilde J_1=\int_0^{\infty} \left[\int_{\frac{-x_1}{2}}^{x_1}\left (
x_1(x_1-x_2)(2x_1+x_2)(x_1+2x_2)e^{-3(x_1^2+x_2^2+x_1x_2)}\right )
dx_2\right ]dx_1=1/18. \]

It follows that when $k=3$, $\tilde
c^{(3,0)}_1=c_{1,E}^{(1;3,0)}=\frac{ 3}{2\sqrt
{\pi}}=0.846284...~~$ Again, $\tilde c^{(3,0)}_2=0$ so $\tilde
c^{(3,0)}_3=-\tilde c^{(3,0)}_1$. Thus, when $k=3$,
\[\tilde \lm^{(3,0)} = \lm_{E}^{(1;3,0)}=\left (\frac{n}{3}+\frac{3}{2\sqrt
{\pi}}\sqrt n\;,\;\frac{n}{3}\;,\;\frac{n}{3}-\frac{3}{2\sqrt
{2\pi}} \sqrt n \right )=~~~~~~~~~~~~~~~~~~~~~~~~~~~~~~~\]

\begin{eqnarray}\label{inv3} =\left (\frac{n}{3}+0.846284\sqrt
n\;,\;\frac{n}{3}\;,\;\frac{n}{3}-0.846284 \sqrt n \right
).\end{eqnarray}

{\bf Note:} for $n=200,\,400,\,600,\,800,\,1000,$ and $1200$,
Remark~\ref{hook1} and 'Mathematica' give the following values of
$\tilde c^{(3,0)}_{1}(n)=c_{1,E}^{(1;3,0)}$:\\
$0.789051,\,0.80486,\,0.812115,\,0.816513,\,0.819547$ and
$0.821802$, as predicted by Equation~(\ref{hookeqn3})

\bigskip

\centerline{\bf Acknowledgement:}

I would like to thank G. Olshanski for some very fruitful
discussions and suggestions.

\end{document}